\documentclass[12pt]{elsarticle}
\usepackage{amsthm,amsmath,amssymb,
mathrsfs,MnSymbol,accents}
\usepackage[all]{xy}

\usepackage[active]{srcltx} 

\DeclareMathOperator{\diag}{diag}
\DeclareMathOperator{\ind}{\mathfrak M}

\renewcommand{\le}{\leqslant}
\renewcommand{\ge}{\geqslant}

\newcommand{\ff}{\mathbb F}
\newcommand{\rr}{\mathbb R}
\newcommand{\cc}{\mathbb C}
\newcommand{\hh}{\mathbb H}
\newcommand{\mc}{\mathcal}
\newcommand{\st}{\upY}
\newcommand{\wt}{\widetilde}

\newcommand{\C}{\mathfrak C}

\newcommand{\ma}[1]{\begin{matrix}
 #1  \end{matrix}}

\newcommand{\mat}[1]{\begin{bmatrix}
 #1  \end{bmatrix}}

\newcommand{\matt}[1]{\left[\begin{smallmatrix}
 #1 \end{smallmatrix}\right]}

\newcommand{\dd}{\lcurvearrowse}
\newcommand{\ud}{\updownarrow}

\newcommand{\cim}{
\begin{picture}(6,6)
\put(3,3){\circle*{3}}
\end{picture}}

\newtheorem{theorem}{Theorem}
\newtheorem{lemma}{Lemma}

\theoremstyle{remark}

\theoremstyle{definition}

\begin{document}

\title{Classification of linear operators satisfying
$(Au,v)=(u,A^rv)$ or $(Au,A^rv)=(u,v)$ on a vector space with indefinite scalar product}

\author[zai]{Victor Senoguchi Borges}
\ead{victorsenoguchi@gmail.com}

\author[zai]{Iryna Kashuba}
\ead{kashuba@ime.usp.br}

\author[serg]{Vladimir V.
Sergeichuk}
\ead{sergeich@imath.kiev.ua} \address[serg]
{Institute of Mathematics, Tereshchenkivska
3, Kiev, Ukraine}

\author[zai]{Eduardo Ventilari Sodr\'{e}}
\ead{eduvsodre@usp.br}

\author[zai]{Andr\'{e} Zaidan}
\ead{andre.zaidan@gmail.com}
\address[zai]{Instituto de Matem\'atica e Estat\'istica, Universidade de S\~ao Paulo, Brasil}

\cortext[cor]{Linear Algebra Appl. 2021,  DOI: 10.1016/j.laa.2020.12.005}

\begin{abstract}
We classify all linear
operators $\mc A:V\to V$ satisfying $(\mc
Au,v)=(u,\mc A^rv)$ and all linear
operators satisfying $(\mc
Au,\mc A^rv)=(u,v)$ with $r=2,3,\dots$ on a
complex,  real, or quaternion vector space with scalar product given by a nonsingular symmetric, skew-symmetric, Hermitian, or skew-Hermitian form.
\end{abstract}


\maketitle

\section{Introduction}

Let $\ff$ be $\cc$, $\rr$, or the skew
field of quaternions $\hh$. Let $V$ be a finite dimensional right
vector space over $\ff$ with scalar product
given by a nonsingular form $\mc F:V\times
V\to\ff$ that is symmetric or
skew-symmetric if $\ff\in\{\cc,\rr\}$, and
Hermitian or skew-Hermitian if $\ff\in\{\cc,\hh\}$. Let $r\in \{1,2,\dots\}$.
A linear operator $\mc
A:V\to V$ is \emph{$r$-selfadjoint} if
\[
\mc F(\mc
Au,v)= \mc F(u,\mc A^rv)\quad \text{for all $u,v\in V$};
\]
$\mc
A:V\to V$ is \emph{$r$-unitary} if it is nonsingular and
\[
\mc F(\mc
Au,\mc A^{r}v)= \mc F(u,v)\quad \text{for all $u,v\in V$.}
\]
The $1$-selfadjoint operators are selfadjoint operators on spaces with indefinite scalar product; their classification
is given in \cite{caa,goh_old,goh,lanc,
mehl,mehl1,rod,ser_izv}.
The $1$-unitary operators are unitary operators on spaces with indefinite scalar product;
their classification
is given in \cite{au,goh_old,goh,lanc,mil,
ser_izv,ser_isom}.

We give
canonical matrices of $r$-selfadjoint
operators and $r$-unitary operators for $r\ge 2$. We use the method developed in  \cite{ser_izv}, which reduces the problem of classifying systems of forms and linear mappings to the problem of classifying systems of linear mappings. This method allows to consider the problems of classifying $r$-selfadjoint
operators and $r$-unitary operators as the same classification problem.

Later on, we use the term ``$(-r)$-selfadjoint operators'' instead of ``$r$-unitary operators'' and solve the problem of classifying $r$-selfadjoint operators for each $r\in\mathbb Z\setminus\{-1,0,1\}$.
In matrix form, this problem is formulated as follows: we consider pairs $(A,F)$ of $n\times n$ matrices
over $\cc$ or $\rr$ satisfying
\begin{equation}\label{666}
A^TF=FA^r,\qquad F^T=F\text{ is nonsingular}
\end{equation}
and give their canonical form
with respect to transformations
\begin{equation}\label{ccdf}
(A,F)\mapsto (S^{-1}AS,S^TFS),\qquad S\text{ is nonsingular};
\end{equation}
we also consider matrix pairs $(A,F)$
over $\cc$ or $\hh$ satisfying
\begin{equation}\label{777}
A^*F=FA^r,\qquad F^*=F\text{ is nonsingular}
\end{equation}
and give their canonical form
with respect to transformations
\begin{equation}\label{nhx}
(A,F)\mapsto (S^{-1}AS,S^*FS),\qquad S\text{ is nonsingular}
\end{equation}
($A$ is nonsingular if $r<0$, and $S^*:=\overline{S}^{\,T}$).

This research was inspired by the articles
\cite{cat1,cat,cat2,LL2,LL3,LL,leb,LL1}, in which Catral, Lebtahi, Romero,
Stuart, Thome, and Weaver study
$\{R,s+1,k\}$-potent (respectively, $\{R,s+1,k,*\}$-potent) matrices; i.e., those matrices $A\in\mathbb C^{n\times n}$ that satisfy $RA=A^{s+1}R$ (respectively, $RA^*=A^{s+1}R$), in which $R\in\mathbb C^{n\times n}$ is a given matrix satisfying $R^k=1$ and $s,k$ are positive integers; compare with \eqref{666} and \eqref{777}.

Each sesquilinear form $\mc F:V\times
V\to \ff$ that we consider
is semilinear in the first argument and linear in the second; $\mc F:V\to
V$ is \emph{skew-Hermitian} if $\mc F(u,v)=-\overline{\mc F(v,u)}$ for all $u,v\in V$. We do not consider skew-Hermitian forms over $\cc$ since if $\mc F(u,v)$ is skew-Hermitian, then $i\mc F(u,v)$ is Hermitian.

Define the matrix
\begin{equation}\label{oky}
(a+bi)^{\rr}:=\mat{a&-b\\b&a}
\qquad\text{for each } a+bi\in\cc \ (a,b\in\rr),
\end{equation}
and the direct sum of matrix pairs
\[
(A_1,F_1)\oplus(A_2,F_2):= \left(\mat{A_1&0\\0&A_2},\mat{F_1&0\\0&F_2}\right).
\]
The notation $\lambda \dd\mu$ means that a parameter $\lambda $ is determined up to replacement by $\mu$. We write $\lambda\in \cc^{\ud}$ if $\lambda$ is a complex parameter that is determined up to replacement by its complex conjugate $\bar\lambda$.
We write ``$(A,\pm F)$'' instead of ``$(A,F)$ and $(A,-F)$''. We denote by $0_n$ and $I_n$ the $n\times n$ zero and identity matrices.

Our main result is the following theorem.

\begin{theorem}\label{ttt}
Let $r\in\mathbb Z\setminus\{-1,0,1\}$.
\begin{itemize}
  \item[\rm(A)]
 Let $V_{\cc}$ be a vector space over $\cc$.
\begin{itemize}
  \item[\rm (a$_1$)] Let
$\mc A$ be an $r$-selfadjoint operator on $V_{\cc}$ with a nonsingular symmetric form $\mc F$. Then
there exists a basis of\/ $V_{\cc}$ in which the pair
$(\mc A,\mc F)$ is given by a direct sum,
uniquely determined up to permutations of
summands, of pairs of the form
\begin{equation*}\label{aa1}
([\lambda ],[1]),\qquad\left(\mat{\mu  &0\\0&\mu
^r},\mat{0&1\\1&0}\right),
\end{equation*}
in which
$\lambda,\mu \in\cc$, $\lambda ^{r}=\lambda$, $\mu^{r^2}=\mu$,
$\mu^{r}\ne \mu$, $\mu\dd\mu^r$.

  \item[\rm (a$_2$)]  Let $\mc A$ be an $r$-selfadjoint
operator on  $V_{\cc}$ with a nonsingular Hermitian form
$\mc F$. Then there exists a basis of\/ $V_{\cc}$ in
which the pair $(\mc A,\mc F)$ is given by
a direct sum, uniquely determined up to
permutations of summands, of pairs of the
form
\begin{equation*}\label{aa3}
([\lambda ],\pm[1]),\qquad\left(\mat{\mu &0\\0&\bar\mu
^r},\mat{0&1\\1&0}\right),
\end{equation*}
in which $\lambda,\mu
\in\cc$, $\lambda ^{r}=\bar \lambda $, $\mu^{r^2}=\mu$,
$\mu^{r}\ne \bar\mu$, $\mu\dd\bar{\mu}^r$.

  \item[\rm (a$_3$)]
Let $\mc A$ be an $r$-selfadjoint
operator on  $V_{\cc}$ with a nonsingular skew-symmetric
form $\mc F$. Then there exists a basis of\/
$V_{\cc}$ in which the pair $(\mc A,\mc F)$ is
given by a direct sum, uniquely determined
up to permutations of summands, of pairs of
the form
\begin{equation*}\label{aa2}
\left(\mat{\lambda&0\\0&\lambda
^r},\mat{0&1\\-1&0}\right),
\end{equation*}
in which $\lambda\in\cc$,
$\lambda^{r^2}=\lambda$, and
$\lambda \dd\lambda ^r$.
\end{itemize}

\item[\rm(B)]
Let $V_{\rr}$ be a vector space over $\rr$.

\begin{itemize}
  \item[\rm (b$_1$)]
Let $\mc A$ be an $r$-selfadjoint
operator on  $V_{\rr}$ with a nonsingular symmetric form
$\mc F$. Then there exists a basis of\/ $V_{\rr}$ in
which the pair $(\mc A,\mc F)$ is given by
a direct sum, uniquely determined up to
permutations of summands, of pairs of the
form
\begin{gather*}
([0],\pm[1]),\quad ([1],\pm[1]),\quad
([-1],\pm[1])\text{ if $r$ is odd},\\
\left(\lambda^{\rr} ,\mat{0&1\\1&0}\right),
          \quad
\left(\mu^{\rr} ,\pm I_2\right),\quad
\left(
\mat{\nu^{\rr}&0\\0&
(\bar{\nu}^r)^{\rr}} ,\mat{0&I_2\\I_2&0}\right),
\end{gather*}
in which $\lambda,\mu,\nu \in\cc^{\ud}\setminus\rr$,
$\lambda ^{r}=\lambda$,
$\mu ^{r}=\bar\mu$,
 $\nu^{r^2}=\nu$, $\nu ^{r}\ne\nu$, $\nu ^{r}\ne\bar\nu$, $\nu \dd\nu ^r$.

 \item[\rm (b$_2$)]
Let $\mc A$ be an $r$-selfadjoint
operator on $V_{\rr}$ with a nonsingular skew-symmetric form
$\mc F$. Then there exists a basis of\/ $V_{\rr}$ in
which the pair $(\mc A,\mc F)$ is given by
a direct sum, uniquely determined up to
permutations of summands, of pairs of the
form
\begin{gather*}
\left(0_2 ,\mat{0&-1\\1&0}\right),\quad
\left(I_2 ,\mat{0&-1\\1&0}\right),\quad
\left(-I_2 ,\mat{0&-1\\1&0}\right)\text{ if $r$ is odd},
             \\
\left(\lambda ^{\rr} ,\mat{0&-1\\1&0}\right),\quad
\quad\left(
\mat{\mu^{\rr}&0\\0&(\bar{\mu}^{r})^{\rr}} ,\mat{0&-I_2\\I_2&0}\right),
\end{gather*}
in which  $\lambda ,\mu
\in\cc^{\ud}\setminus\rr,$ $\lambda ^{r}=\bar \lambda $, $\mu^{r^2}=\mu$, $\mu^{r}\ne \bar\mu$, $\mu\dd\mu^r$.
\end{itemize}

\item[\rm(C)]
Let $V_{\hh}$ be a right vector space over $\hh$.

\begin{itemize}
  \item[\rm (c$_1$)]
Let $\mc A$ be an $r$-selfadjoint
operator on  $V_{\hh}$ with a nonsingular Hermitian form
$\mc F$ with respect to quaternion
conjugation \begin{equation}\label{yyb}
h=a+bi+cj+dk\ \mapsto\ \bar
h=a-bi-cj-dk,\qquad a,b,c,d\in\rr.
\end{equation}
Then there exists a basis of\/ $V_{\hh}$ in
which the pair $(\mc A,\mc F)$ is given by
a direct sum, uniquely determined up to
permutations of summands, of pairs of the
form
\[
([\lambda],\pm [1]),\qquad
\left(\mat{\mu &0\\0&\bar \mu^r} ,\mat{0&1\\1&0}\right),
\]
in which
$\lambda,\mu \in\cc^{\ud},$
$\lambda^{r}=\bar\lambda$,
$\mu^{r^2}=\mu,$
$\mu^{r}\ne \bar\mu,$ $\mu
\dd\mu ^r$.

  \item[\rm (c$_2$)]
Let $\mc A$ be an $r$-selfadjoint
operator on  $V_{\hh}$ with a nonsingular Hermitian form
$\mc F$ with respect to quaternion
semiconjugation \begin{equation}\label{ybb}
h=a+bi+cj+dk\ \mapsto\ \widehat h= a-bi+cj+dk,\qquad
a,b,c,d\in\rr. \end{equation}
Then there exists a basis of\/ $V_{\hh}$ in
which the pair $(\mc A,\mc F)$ is given by
a direct sum, uniquely determined up to
permutations of summands, of pairs of the
form
\begin{gather*}
([\lambda],\pm [1])\text{ if }\lambda\notin\rr,\qquad
([\lambda],[1])\text{ if }\lambda\in\rr,
                            \\
([\mu],[j]),\qquad
\left(\mat{\nu &0\\0&\bar \nu^r} ,\mat{0&1\\1&0}\right),
\end{gather*}
in which
$\lambda,\mu,\nu\in\cc^{\ud},$
$\lambda^{r}=\bar\lambda$,
${\mu}^{r}=\mu\notin\rr$, $\nu ^{r^2}=\nu,$
${\nu}^{r}\ne \nu$,
${\nu}^{r}\ne \bar\nu$,
$\nu \dd\nu^r$.

  \item[\rm (c$_3$)]
Let $\mc A$ be an $r$-selfadjoint
operator on  $V_{\hh}$ with a nonsingular form
$\mc F$ that is skew-Hermitian with respect to quaternion
conjugation \eqref{yyb}. Then there exists a basis of\/ $V_{\hh}$ in
which the pair $(\mc A,\mc F)$ is given by
a direct sum, uniquely determined up to
permutations of summands, of pairs of the
form
\begin{gather*}
([\lambda],\pm [i])\text{ if }\lambda\notin\rr,\qquad
([\lambda],[i])\text{ if }\lambda\in\rr,
                            \\
([\mu],[j]),\qquad
\left(\mat{\nu &0\\0&\bar \nu^r} ,\mat{0&-1\\1&0}\right),
\end{gather*}
in which
$\lambda,\mu,\nu\in\cc^{\ud}$,
$\lambda^{r}=\bar\lambda$,
$\mu\notin\rr,$
${\mu}^{r}=\mu$, $\nu ^{r^2}=\nu$,
${\nu}^{r}\ne \nu$,
${\nu}^{r}\ne \bar\nu$,
$\nu
\dd\nu^r$.

  \item[\rm (c$_4$)]
Let $\mc A$ be an $r$-selfadjoint
operator on  $V_{\hh}$ with a nonsingular form
$\mc F$ that is skew-Hermitian with respect to quaternion
semiconjugation  \eqref{ybb}. Then there exists a basis of\/ $V_{\hh}$ in
which the pair $(\mc A,\mc F)$ is given by
a direct sum, uniquely determined up to
permutations of summands, of pairs of the
form
\[
([\lambda],\pm [i]),\qquad
\left(\mat{\mu &0\\0&\bar \mu^r} ,\mat{0&-1\\1&0}\right),
\]
in which
$\lambda,\mu\in\cc^{\ud}$,
$\lambda^{r}=\bar\lambda$, $\mu ^{r^2}=\mu$,
$\mu^{r}\ne \bar\mu$, $\mu
\dd\mu ^r$.
\end{itemize}
\end{itemize}
\end{theorem}

Each condition $\lambda ^{r^2}=\lambda $,
$\lambda ^{r}=\lambda $, or
$\lambda ^{r}=\bar\lambda $ implies that $\lambda \ne 0$ if $r<0$. Theorem \ref{ttt} remains true if $\cc$, $\rr$, and $\hh$ are replaced by an algebraically closed field of zero characteristic, a real closed field, and the skew field of quaternions over a real closed field, respectively.

An \emph{involution} $a\mapsto\wt a$ on a field or skew field $\ff$ is a bijection $\ff\to\ff$
satisfying \[ \wt{a+b}=\wt a+\wt b,\quad
\wt{\,ab\,}=\wt b\ \wt a,\quad \wt{\wt
a}=a\qquad\text{for all }a,b\in\ff. \]
If an involution on
$\hh$ is not quaternion
conjugation \eqref{yyb}, then it is quaternion
semiconjugation \eqref{ybb} in a suitable set of the fundamental units
$i,j,k$; see \cite[Lemma 2.2]{ser_isom}.

\section{Reduction of the problem of classifying $r$-selfadjoint
operators to the problem of classifying matrices under similarity}

We prove
Theorem  \ref{ttt} in the next section by the method that is developed in
\cite{ser_izv}. It
reduces the
problem of classifying systems of linear
mappings and forms to the problem of
classifying systems of linear mappings. Bilinear and sesquilinear forms, pairs of
symmetric, skew-symmetric, and Hermitian
forms, unitary and selfadjoint operators on
a vector space with indefinite scalar product
are classified in \cite{ser_izv} over a field $\mathbb K$ of characteristic not 2 up to
classification of Hermitian forms over
finite extensions of $\mathbb K$ (and so they are fully classified over $\rr$ and
$\cc$).

The reader is expected to be familiar with this method; it is described in details in  \cite{ser_isom} and is used
in \cite{hor-ser_can1,
hor-ser,mel,ser_sur}. In this section,
we sketchily describe it in a special case: for the problem of classifying $r$-selfadjoint
operators.

Systems consisting of vector spaces and of linear
mappings and forms on them are considered as \emph{representations of mixed graphs}; i.e., graphs with undirected and directed edges.
Its vertices represent vector spaces, its
undirected edges represent forms, and its directed edges represent linear mappings.

In particular, each pair $(\mc A,\mc F)$ from
Theorem \ref{ttt} defines the representation
\begin{equation}\label{rty1}
\xymatrix{
 {V\,}
\save !<-6pt,0pt>
 \ar@(ul,dl)@{<-}_{\mc A}
\restore \save !<6pt,0pt>
 \ar@<-0.4ex>@(ur,dr)@{-}^{\mc F}
\restore}\qquad \begin{matrix} \mc F(\mc
Au,v)= \mc F(u,\mc A^rv)
\text{ if } r\ge 2,\\
\mc F(\mc
Au,\mc A^{-r}v)= \mc F(u,v)
\text{ if } r\le -2,\\
\mc F(u,v)=
\varepsilon \wt{\mc F(v,u)}\text{ is
nonsingular}
  \end{matrix}
\end{equation}
of the mixed graph
$
\xymatrix{
 {\bullet}
\save !<-2pt,0pt>
 \ar@(ul,dl)@{<-}_{}
\restore \save !<2pt,0pt>
 \ar@<-0.4ex>@(ur,dr)@{-}^{}
\restore}
$
over $\ff$ with involution $a\mapsto\wt a$, in which $\varepsilon:=1$ if
$\mc F$ is symmetric or Hermitian, and
$\varepsilon:=-1$ if $\mc F$ is
skew-symmetric or skew-Hermitian. Choosing a basis in $V$, we
give \eqref{rty1} by its matrices
\begin{equation}\label{rty2}
\xymatrix{
 {n\,}
\save !<-5pt,0pt>
 \ar@(ul,dl)@{<-}_{A}
\restore \save !<5pt,0pt>
 \ar@<-0.4ex>@(ur,dr)@{-}^{F}
\restore}\qquad A^{\st}F=FA^r,\quad
F^{\st}=\varepsilon  F\text{ is
nonsingular,} \end{equation} in which
$n:=\dim V$, $A$ and $F$ are $n\times n$
matrices  over $\ff$, and $A^{\st}:=\wt
A^T$ ($A^{\st}=A^T$ if $a\mapsto\wt a$ is the identity  involution, and $A^{\st}=A^*$ otherwise). Changing the basis in $V$, we can
reduce $(A,F)$ by transformations
\begin{equation}\label{kix}
(A,F)\mapsto(S^{-1}AS,S^{\st}FS),\qquad
S\text{ is nonsingular}
\end{equation}
(see \eqref{ccdf} and \eqref{nhx}).
We
say that the pairs $(A,F)$ and
$(S^{-1}AS,S^{\st}FS)$ are \emph{isomorphic
via $S$}.

Replacing $\mc F:V\times V\to \ff$ in
\eqref{rty1} by the pair of mutually
adjoint linear mappings $\mc F:v\mapsto \mc
F(?,v)$ and $\mc F^{\st}:u\mapsto \widetilde{\mc
F(u,?)}$, we obtain the system of linear
mappings
\begin{equation}\label{qwe}
\xymatrix{
 {V}
\save !<-4pt,0pt>
 \ar@(ul,dl)@{<-}_{\mc A}
\restore
  \ar@/^/@{->}[rr]^{\mc F}
 \ar@/^/@{<-}[rr];[]^{\mc F^{\st}}
  &&{V^{\st}}
\save !<8pt,0pt> \ar@(ur,dr)@{->}^{\mc
A^{\st}} \restore
 }\quad
\mc A^{\st}\mc F=\mc F\mc A^r,\ \ \mc F^{\st}=\varepsilon
\mc F\text{ is nonsingular,}
\end{equation}
in which $V^{\st}$ is the \emph{$^{\st\!}$dual space} (with respect to the
involution $a\mapsto\wt a$) consisting of
semilinear forms  on $V$, and
$\mc A^{\st}:V^{\st}\to V^{\st}$ is the \emph{$^{\st\!}$dual mapping} defined by $\varphi\mapsto \varphi\mc A$. In the
matrix form, \begin{equation}\label{qwe1}
\xymatrix{
 {n}
\save !<-4pt,0pt>
 \ar@(ul,dl)@{<-}_{A}
\restore
  \ar@/^/@{->}[rr]^{F}
 \ar@/^/@{<-}[rr];[]^{F^{\st}}
  &&{n}
\save !<4pt,0pt> \ar@(ur,dr)@{->}^{A^{\st}}
\restore
 }\quad
A^{\st}F=FA^r,\ \ F^{\st}=\varepsilon
F\text{ is nonsingular.}
\end{equation}

Thus, there is the bijective correspondence
\begin{equation}\label{rvr}
\xymatrix{
 {n\,}
\save !<-5pt,0pt>
 \ar@(ul,dl)@{<-}_{A}
\restore \save !<5pt,0pt>
 \ar@<-0.4ex>@(ur,dr)@{-}^{F}
\restore}\qquad \mapsto\qquad
\xymatrix{
 {n}
\save !<-4pt,0pt>
 \ar@(ul,dl)@{<-}_{A}
\restore
  \ar@/^/@{->}[rr]^{F}
 \ar@/^/@{<-}[rr];[]^{F^{\st}}
  &&{n}
\save !<4pt,0pt> \ar@(ur,dr)@{->}^{A^{\st}}
\restore
 }
\end{equation}
between the matrix sets of systems \eqref{rty1} and
\eqref{qwe}.

Let us consider a system of linear mappings over
$\ff$:
\begin{equation}\label{ndu}
\mc M:\ \xymatrix{
 {V_1}
\save !<-2mm,0cm>
 \ar@(ul,dl)@{<-}_{\mc A_1}
\restore
  \ar@/^/@{->}[rr]^{\mc F_1}
 \ar@/^/@{<-}[rr];[]^{\mc F_2}
  &&{V_2}
\save !<2mm,0cm> \ar@(ur,dr)@{->}^{\mc A_2}
\restore
 }\qquad
 {\begin{matrix}
\mc A_2\mc F_1=\mc F_1\mc A_1^r, \ \mc A_2^r\mc F_2=\mc F_2\mc A_1,\\
\mc  F_1=\varepsilon\mc  F_2\text{ is nonsingular},
 \end{matrix}}
\end{equation}
which is a representation of the quiver
$\xymatrix@=15pt{
&*{\cim}
 \ar@(ul,dl)@{<-}
  \ar@/^/@{->}[r]
 \ar@/^/@{<-}[r];[]
  &*{\cim}
\ar@(ur,dr)@{->}
}$\ \quad\ . Choosing bases in $V_1$ and $V_2$, we give it by
a system of $n\times n$ matrices ($n:=\dim V_1=\dim V_2$)
\begin{equation}\label{sju}
M:\ \xymatrix{
 {n}
\save !<-1mm,0cm>
 \ar@(ul,dl)@{<-}_{A_1}
\restore
  \ar@/^/@{->}[rr]^{F_1}
 \ar@/^/@{<-}[rr];[]^{F_2}
  &&{n}
\save !<1mm,0cm> \ar@(ur,dr)@{->}^{A_2}
\restore
 }\qquad
 {\begin{matrix}
A_2F_1=F_1A_1^r, \ A_2^rF_2=F_2A_1,\\
 F_1=\varepsilon F_2\text{ is nonsingular}.
 \end{matrix}}
\end{equation}
Changing bases in $V_1$ and $V_2$, we can reduce \eqref{sju} by transformations
\begin{equation}\label{smk}
M':\ \xymatrix{
 {n}
\save !<-1mm,0cm>
 \ar@(ul,dl)@{<-}_{RA_1R^{-1}}
\restore
  \ar@/^/@{->}[rr]^{SF_1R^{-1}}
 \ar@/^/@{<-}[rr];[]^{SF_2R^{-1}}
  &&{n}
\save !<1mm,0cm>
\ar@(ur,dr)@{->}^{SA_2S^{-1}} \restore
 }\qquad
R,S\text{ are nonsingular.}
 \end{equation}
We
say that the matrix sets \eqref{sju} and
\eqref{smk} are \emph{isomorphic via $R$ and $S$} and write
$M\simeq M'$.
This isomorphism can be shown by the commutative diagram
\begin{equation*}\label{mcm}
\xymatrix@C=50pt{
{n} \ar[r]^{A_1}\ar[d]_{R}
&{n}  \ar[d]_{R}
\ar@<-0.4ex>[r]_{F_2}
 \ar@<0.4ex>[r]^{F_1}
&{n}     \ar[d]^{S}
 \ar[r]^{A_2}
&{n}\ar[d]^{S}
   \\
{n} \ar[r]^{RA_1R^{-1}}
&{n}
\ar@<-0.4ex>[r]_{SF_2R^{-1}}
 \ar@<0.4ex>[r]^{SF_1R^{-1}}
&{n}
 \ar[r]^{SA_2S^{-1}}
&{n}
 }\end{equation*}

The \emph{direct sum} of systems $M$ and
$M'$ is the system
\begin{equation}\label{mju}
M\oplus
M':\quad\xymatrix{
 {n+n'}
\save !<-5mm,0cm>
 \ar@(ul,dl)@{<-}_{A_1\oplus A_1'}
\restore
  \ar@/^/@{->}[rr]^{F_1\oplus F_1'}
 \ar@/^/@{<-}[rr];[]^{F_2\oplus F_2'}
  &&{n+n'}
\save !<5mm,0cm>
\ar@(ur,dr)@{->}^{A_2\oplus A_2'} \restore
 }
\end{equation}
A system is \emph{indecomposable} if it is not isomorphic to a direct sum of the form \eqref{mju} with nonzero $n$ and $n'$.

For each system \eqref{sju},
we define the \emph{dual system}
\begin{equation*}\label{sju1}
M^{\circ}:\quad\xymatrix{
 {n}
\save !<-1mm,0cm>
 \ar@(ul,dl)@{<-}_{A_2^{\st}}
\restore
  \ar@/^/@{->}[rr]^{F_2^{\st}}
 \ar@/^/@{<-}[rr];[]^{F_1^{\st}}
  &&{n}
\save !<1mm,0cm>
\ar@(ur,dr)@{->}^{A_1^{\st}} \restore
 }
\end{equation*}
A system $M$ is
\emph{selfdual} if $M=M^{\circ}$, which means that it
has the form \eqref{qwe1}.

Suppose we know the following sets of systems of the form \eqref{sju}:
\begin{description}
  \item[$\ind_{\varepsilon}(\ff)$] which is a set of nonisomorphic indecomposable systems such that every
indecomposable system \eqref{sju} is
isomorphic to exactly one system from
$\ind_{\varepsilon}(\ff)$,

  \item[$\ind'_{\varepsilon}(\ff)$]
which is a set of nonisomorphic indecomposable
selfdual systems such that every
indecomposable selfdual system is
isomorphic to exactly one system from
$\ind'_{\varepsilon}(\ff)$,

  \item[$\ind''_{\varepsilon}(\ff)$]
which is a set of nonisomorphic indecomposable systems that are not isomorphic to selfdual such that every
indecomposable system that is not isomorphic to selfdual is
isomorphic to exactly one system from
$\ind''_{\varepsilon}(\ff)$.
\end{description}

For each
$M\in\ind'_{\varepsilon}(\ff)$ of the form \eqref{qwe1}, we
define the matrix pairs
\begin{equation}\label{re2}
\mathring
M:\ \xymatrix{
 {n\,}
\save !<-3pt,0pt>
 \ar@(ul,dl)@{<-}_{A}
\restore \save !<3pt,0pt>
 \ar@<-0.4ex>@(ur,dr)@{-}
 ^{F}
\restore}\qquad\qquad \mathring
M^-:\ \xymatrix{
 {n\,}
\save !<-3pt,0pt>
 \ar@(ul,dl)@{<-}_{A}
\restore \save !<3pt,0pt>
 \ar@<-0.4ex>@(ur,dr)@{-}
 ^{-F}
\restore} \end{equation}
and for each $N\in\ind''_{\varepsilon}(\ff)$ of the form
\eqref{sju}, we define the matrix pair
\begin{equation*}\label{rk2}
N^+:\quad\xymatrix{
 {2n\,}
\save !<-8pt,0pt>
 \ar@(ul,dl)@{<-}_{\mat{A_1&0\\0&A_2^{\st}}}
\restore \save !<8pt,0pt>
 \ar@<-0.4ex>@(ur,dr)@{-}
 ^{\mat{0&F_2^{\st}\\F_1&0}}
\restore}.
\end{equation*}
Note that the natural bijection \eqref{rvr} takes $\mathring M$ into $M$, and $N^+$ into a selfdual system that is isomorphic to $N\oplus N^{\circ}$.

The following lemma reduces the problem of classifying matrix pairs \eqref{rty2} up to transformations \eqref{kix} to the problem of classifying systems \eqref{sju} up to transformations \eqref{smk}.
This lemma is a special case of
\cite[Theorem 3.2]{ser_isom} about
arbitrary systems of linear mappings and
forms.

\begin{lemma}\label{ll} Each pair
\eqref{rty2} over $\ff\in\{\cc,\rr,\hh\}$
is isomorphic to a direct sum of
pairs of the types
\[
  N^+\quad  \text{and}\quad
\begin{cases}
   \mathring{M} &\text{if $\mathring{M}^{-}$ and
   $\mathring{M}$ are isomorphic}, \\
\mathring{M},\ \mathring{M}^{-} &\text{if $\mathring{M}^{-}$ and
   $\mathring{M}$ are not isomorphic},
  \end{cases}
\]
in which  $M\in\ind'_{\varepsilon}(\ff)$
and $N\in\ind''_{\varepsilon}(\ff)$. This sum is uniquely determined,
up to permutations of direct summands and replacements of $N\in\ind''_{\varepsilon}(\ff)$ by $N^{\circ}$.
\hfill\hbox{\qedsymbol} \end{lemma}

\begin{lemma}\label{lll}
If a system of the form \eqref{sju} is isomorphic to a selfdual system via $R$ and $S$, then it is isomorphic to some selfdual system via $I$ and $R^{\st}S$.
\end{lemma}

\begin{proof} The corresponding selfdual system is constructed as follows:
\[ \xymatrix{
 {n}  \ar[d]_{R}
\save !<-1pt,0pt>
 \ar@(ul,dl)@{<-}_{A_1}
\restore \ar@<-0.4ex>[rr]_{F_2}
 \ar@<0.4ex>[rr]^{F_1}
  &&{n}     \ar[d]^{S}
\save !<1pt,0pt>
\ar@(ur,dr)@{->}^{A_2}
\restore
   \\
 {n} \ar[d]_{R^{-1}}
\save !<-1pt,0pt>
 \ar@(ul,dl)@{<-}_{B}
\restore \ar@<-0.4ex>[rr]_{G^{\st}}
 \ar@<0.4ex>[rr]^{G}
  &&{n} \ar[d]^{R^{\st}}
\save !<1pt,0pt>
\ar@(ur,dr)@{->}^{B^{\st}
\ \text{ - selfdual}}
\restore
   \\
 {n}
\save !<-1pt,0pt>
 \ar@(ul,dl)@{<-}_{R^{-1}BR}
\restore \ar@<-0.4ex>[rr]_{R^{\st}G^{\st}R}
 \ar@<0.4ex>[rr]^{R^{\st}GR}
  &&{n}
\save !<1pt,0pt>
\ar@(ur,dr)@{->}^{R^{\st}B^{\st}(R^{\st})^{-1}
\ \text{ - selfdual}}
\restore
 }
\]
\vskip-2em
\end{proof}

\section{Proof of Theorem \ref{ttt}}

\begin{lemma}\label{dcj}
Each square matrix over $\cc$, $\rr$, and $\hh$ is similar to a direct sum, uniquely determined up to permutations of summands, of matrices from the following matrix sets:
\begin{itemize}
  \item[\rm(a)]
$ \mathfrak
C(\cc):=\{J_n(\lambda
)\,|\,\lambda\in\cc\},$ in which
\begin{equation*}\label{888} J_n(\lambda
):=\mat{\lambda &1&&0\\ &\lambda
&\ddots\\&&\ddots&1\\0&&&\lambda }\
(n\text{-by-}n).
\end{equation*}

  \item[\rm(b)]
$
\C(\rr)=:\{J_n(a)\,|\,a\in\rr\}\cup\{J_n(\lambda
)^{\rr}\,|\,\lambda
\in\cc^{\ud}\setminus\rr\},
$ in which
\[J_n(\lambda
)^{\rr}:=\mat{\lambda^{\rr} &I_2&&0\\
&\lambda^{\rr}
&\ddots\\&&\ddots&I_2\\0&&&\lambda^{\rr} }\
(2n\text{-by-}2n),\]
$\lambda^{\rr}$ is determined in \eqref{oky}, and $
\cc^{\ud}=\{\lambda \in\cc\,|\,\lambda \dd\bar\lambda\}$.

  \item[\rm(c)]
$
\C(\hh):=\{J_n(\lambda)\,|\,\lambda
\in\cc^{\ud}\}.
$
\end{itemize}
\end{lemma}

\begin{proof} The statement (a) is the Jordan theorem;
(b) and (c) are given in \cite[Theorem 3.4.1.5]{hor-joh} and \cite[Theorem 5.5.3]{rod}.
\end{proof}

Each system \eqref{sju} is reduced by
transformations \eqref{smk} with $R=F_1$ and
$S=I_n$ to a system of the form
\begin{equation*}\label{qook} M_{\varepsilon
}(A):\quad\xymatrix{
 {n}
\save !<-2pt,0pt>
 \ar@(ul,dl)@{<-}_{A}
\restore
  \ar@/^/@{->}[rr]^{I_n}
 \ar@/^/@{<-}[rr];[]^{\varepsilon I_n}
  &&{n}
\save !<2pt,0pt> \ar@(ur,dr)@{->}^{A^r}
\restore
 }\qquad\qquad
A^{r^2}=A,
 \end{equation*}
whose dual system is
\begin{equation*}\label{qooy} M_{\varepsilon
}(A)^{\circ}:\quad\xymatrix{
 {n}
\save !<-2pt,0pt>
 \ar@(ul,dl)@{<-}_{(A^r)^{\st}}
\restore
  \ar@/^/@{->}[rr]^{\varepsilon I_n}
 \ar@/^/@{<-}[rr];[]^{I_n}
  &&{n}
\save !<2pt,0pt> \ar@(ur,dr)@{->}^{A^{\st}}
\restore
 }
 \end{equation*}
Clearly, $M_{\varepsilon }(A)\simeq
M_{\varepsilon }(B)$ if and only if $A$ and
$B$ are similar; and so
\begin{equation}\label{tct}
\ind_{\varepsilon}(\ff)=\{M_{\varepsilon
}(A)\,|\,A\in \mathfrak  C( \ff)\text{ such
that }A^{r^2}=A\};
\end{equation}
\begin{equation}\label{fbp}
M_{\varepsilon }(B)\simeq
M_{\varepsilon
}(A)^{\circ}\quad\Longleftrightarrow\quad
\text{$B$ is similar to $(A^r)^{\st}$}.
\end{equation}

\subsection{Case {\rm(A)}: $\ff=\cc$}

Let $a\mapsto \wt a$ be the identity involution or complex conjugation.
Suppose that
$J_n(\lambda)^{r^2}=J_n(\lambda)$ with $\lambda \in\cc$. If
$\lambda =0$, then $n=1$ since $r^2\ge 4$. If
$\lambda \ne 0$, then
$J_n(\lambda)^{r^2-1}=I_n$, all entries of
the first over-diagonal of
$J_n(\lambda)^{r^2-1}$ are $(r^2-1)\lambda
^{r^2-2}$, and so $n=1$ too. Thus,
\begin{equation}\label{oof}
J_n(\lambda)^{r^2}=J_n(\lambda)\quad\Longrightarrow
\quad n=1
\end{equation}
and
\[\ind_{\varepsilon}(\cc)=\{M_{\varepsilon
}(\lambda)\,|\,\lambda \in\cc,\
\lambda^{r^2}=\lambda \}\] (to simplify
notation, we write $\lambda$ instead of
$\lambda I_1$).

Since \begin{equation*}\label{kmj}
M_{\varepsilon }(\lambda):\: \xymatrix{
 {1}
\save !<-2pt,0pt>
 \ar@(ul,dl)@{<-}_{\lambda}
\restore
  \ar@/^/@{->}[r]^{1}
 \ar@/^/@{<-}[r];[]^{\varepsilon}
  &{1}
\save !<2pt,0pt>
\ar@(ur,dr)@{->}^{\lambda^r} \restore
 }\qquad
 M_{\varepsilon }(\lambda)^{\circ}:\:
 \xymatrix{
 {1}
\save !<-2pt,0pt>
 \ar@(ul,dl)@{<-}_{\wt{\lambda}^r}
\restore
  \ar@/^/@{->}[r]^{\varepsilon}
 \ar@/^/@{<-}[r];[]^{1}
  &{1}
\save !<2pt,0pt>
\ar@(ur,dr)@{->}^{\wt{\lambda}} \restore
 }
\end{equation*}
and $\lambda
^{r^2}=\lambda$, we have $
M_{\varepsilon }(\lambda)^{\circ}\simeq
M_{\varepsilon }(\wt{\lambda}^r)$. Hence,
\begin{equation}\label{fdep}
M_{\varepsilon }(\mu)\simeq
M_{\varepsilon
}(\lambda)^{\circ}\quad\Longleftrightarrow\quad
\mu =\wt{\lambda}^{\;r}.
\end{equation}

\noindent
The following cases are possible:

\begin{itemize}
  \item[(a$_{1}$):] \emph{$\varepsilon =1$ and the involution on $\cc$ is the identity.} Then
\[
M_1(\lambda)\simeq M_1(\lambda)^{\circ}
\quad\Longleftrightarrow\quad M_1(\lambda)=
M_1(\lambda)^{\circ}\quad
\Longleftrightarrow\quad\lambda
={\lambda}^{r},
\]
and so
\begin{align*}
&\ind_1'(\cc)=\{ M_1(\lambda)\,|\,\lambda\in\cc,\
\lambda ^r=\lambda\,\},
         \\
&\ind_1''(\cc)=\{M_1(\mu)\,|\,\mu \in\cc,\
\mu^{r^2}=\mu,\ \mu^r\ne \mu\}.
\end{align*}
Lemma \ref{ll} and \eqref{fdep} ensure (a$_1$).

  \item[(a$_{2}$):] \emph{$\varepsilon =1$ and the involution on $\cc$ is complex conjugation.} Then
\[
M_1(\lambda)\simeq M_1(\lambda)^{\circ}
\quad\Longleftrightarrow\quad M_1(\lambda)=
M_1(\lambda)^{\circ}\quad
\Longleftrightarrow\quad\lambda
=\bar{\lambda}^{r},
\]
and so
\begin{align*}
&\ind_1'(\cc)=\{ M_1(\lambda)\,|\,\lambda\in\cc,\
\lambda ^r=\bar{\lambda}\},
           \\
&\ind_1''(\cc)=\{M_1(\mu)\,|\,\mu \in\cc,\
\mu^{r^2}=\mu,\ \mu^r\ne \bar{\mu}\}.
\end{align*}
Lemma \ref{ll} and \eqref{fdep} ensure (a$_2$).

  \item[(a$_3$):]  \emph{$\varepsilon =-1$ and the involution on $\cc$ is the identity.}
The system $M_{-1}(\lambda )$ is not
isomorphic to a selfdual system since there
are no nonsingular $1\times 1$ matrices $R$
and $S$ such that
$SI_1R^{-1}=S(-I_1)R^{-1}$ (see
\eqref{smk}). Therefore,
\begin{align*}
\ind_{-1}'(\cc)&=\emptyset,\\
\ind_{-1}''(\cc)&=\{M_{-1}(\lambda)
\,|\,
\lambda \in\cc,\
\lambda^{r^2}=\lambda\}.
\end{align*}
Lemma \ref{ll} and \eqref{fdep} ensure (a$_3$).
\end{itemize}

\subsection{Case {\rm(B)}: $\ff=\rr$}
The set $\C(\rr)$ is given in Lemma \ref{dcj}(b).
The equality \begin{equation}\label{kas}
\begin{split}
\ind_{\varepsilon }(\rr)=&\{M_{\varepsilon
}(0),\ M_{\varepsilon }(1)\} \cup
\{M_{\varepsilon }(-1)\,|\,
\text{if $r$ is odd}\}\\
    & \cup
\{M_{\varepsilon }(\lambda^{\rr})
\,|\,\lambda \in\cc^{\ud}\setminus\rr,\
\lambda^{r^2}=\lambda \}
\end{split}
\end{equation}
is proved as follows:
\begin{itemize}
  \item
Consider $J_n(a)\in \C(\rr)$ with $a\in\rr$ and
$J_n(a)^{r^2}=J_n(a)$. By \eqref{oof},
$n=1$. Since $a$ is real, $a^{r^2}=a$ implies that either $a=0$, or $a=\pm 1$ if $r$ is odd and $a=1$ if $r$ is even.
Note that
\begin{equation}
\label{vfo}
\parbox[c]{0.84\textwidth}{each system $M_{\varepsilon}(\lambda )$
with $\lambda \in\{0,\; 1,\; -1\; (\text{if $r$ is odd})\}$
is selfdual if $\varepsilon =1$; it is not isomorphic
to selfdual if $\varepsilon =-1$.}
\end{equation}

  \item Consider $J_n(\lambda )^{\rr}\in\C(\rr)$ with  $\lambda
\in\cc\setminus\rr$ and $(J_n(\lambda
)^{\rr})^{r^2}=J_n(\lambda )^{\rr}$. The
matrix $J_n(\lambda )^{\rr}$ is similar
over $\cc$ to $J_n(\lambda )\oplus
J_n(\bar{\lambda })$. Hence $J_n(\lambda
)^{r^2}=J_n(\lambda )$, $n=1$,
and $\lambda^{r^2}=\lambda $.
\end{itemize}

For every $M_{\varepsilon
}(\lambda^{\rr}),
M_{\varepsilon }(\mu ^{\rr})\in\ind_{\varepsilon }(\rr)$,
we have
\begin{equation}\label{bxo}
M_{\varepsilon }(\mu ^{\rr})
\simeq
M_{\varepsilon
}(\lambda^{\rr})^{\circ}
\qquad
\Longleftrightarrow\qquad \mu =\lambda^r\text{ or }\mu =\bar{\lambda}^r
\end{equation}
since
 $M_{\varepsilon }(\mu ^{\rr})
\simeq
M_{\varepsilon
}(\lambda^{\rr})^{\circ}
$ if and only if $\mu^{\rr}$ is similar to $((\lambda
^{\rr})^r)^T$, if and only if $\diag(\mu,\bar{\mu})$
is similar to $\diag(\lambda,\bar{\lambda})^r$, if and only if
$\mu=\lambda^r $ or $\mu=\bar{\lambda}^r$.

Thus,
if $M_{\varepsilon }(\lambda^{\rr})$ is
isomorphic to a selfdual system, then
$\lambda^r=\lambda$ or
$\lambda^r=\bar{\lambda}$.
Write $\lambda=a+bi$ $(a,b\in\rr,\ b\ne 0)$, then
$\lambda^{\rr}=\matt{a&-b\\b&a}$ and
\begin{align*}\label{aw1}
\lambda^r=\lambda\quad&\Longrightarrow\quad
\matt{a&-b\\b&a}^r=\matt{a&-b\\b&a},\\
\lambda^r=\bar{\lambda}\quad&\Longrightarrow\quad
\matt{a&-b\\b&a}^r=\matt{a&b\\-b&a}.
\end{align*}

\noindent
The following two cases are possible:
\begin{itemize}
  \item[(b$_1$):]
\emph{$\varepsilon =1$.} Let $\lambda=a+bi$ $(a,b\in\rr,\ b\ne 0)$. If $M_1(\lambda^{\rr})$ is
isomorphic to a selfdual system, then
$\lambda^r=\lambda$ or
$\lambda^r=\bar{\lambda}$. If
$\lambda^r=\lambda$, then
$\matt{a&-b\\b&a}Z =Z\matt{a&b\\-b&a}$ with
$Z:=\matt{0&1\\1&0}$, and so
$M_1(\lambda^{\rr})$ is isomorphic to a
selfdual system:
\begin{equation*}\label{s7u1} \xymatrix{
 {2}  \ar[d]_{I_2}
\save !<-2pt,0pt>
 \ar@(ul,dl)@{<-}_{\matt{a&-b\\b&a}}
\restore \ar@<-0.4ex>[rr]_{I_2}
 \ar@<0.4ex>[rr]^{I_2}
  &&{2}     \ar[d]^{Z}
\save !<2pt,0pt>
\ar@(ur,dr)@{->}^{\matt{a&-b\\b&a}}
\restore
   \\
 {2}
\save !<-2pt,0pt>
 \ar@(ul,dl)@{<-}_{\matt{a&-b\\b&a}}
\restore \ar@<-0.4ex>[rr]_{Z}
 \ar@<0.4ex>[rr]^{Z}
  &&{2}
\save !<2pt,0pt>
\ar@(ur,dr)@{->}^{\matt{a&b\\-b&a}}
\restore
 }
\end{equation*}
If $\lambda^r=\bar{\lambda}$, then the system
$M_1(\lambda^{\rr})$ is selfdual.

Thus,
\begin{align*}\label{yua}
\ind'_1(\rr)=&\{M_1(0),\ M_1(1)\}\cup
\{M_{1}(-1)\,|\,
\text{if $r$ is odd}\}\\
&\cup\big\{\xymatrix{
 {2}
\save !<-2pt,0pt>
 \ar@(ul,dl)@{<-}_{\ma{\lambda^{\rr}}}
\restore \ar@<-0.4ex>[r]_{Z}
 \ar@<0.4ex>[r]^{Z}
  &{2}
\save !<2pt,0pt>
\ar@(ur,dr)@{->}^{\ma{\bar{\lambda}^{\rr}}}
\restore}
\big|\,
\lambda \in\cc^{\ud}\setminus\rr,\
\lambda ^{r}=\lambda
 \big\}\\
&
\cup
\{M_1(\mu^{\rr})\,|\,\mu
\in\cc^{\ud}\setminus\rr,\ \mu^{r}=\bar \mu \},\\
\ind''_1(\rr)
=&\{M_1(\nu^{\rr})\,|\,
\nu\in\cc^{\ud}\setminus\rr,\ \nu
^{r^2}=\nu,\ \nu ^{r}\ne\nu,\  \nu ^{r}\ne\bar\nu\}.
\end{align*}
Each system
$\xymatrix{
 {2}
\save !<-2pt,0pt>
 \ar@(ul,dl)@{<-}_{\ma{\lambda^{\rr}}}
\restore \ar@<-0.4ex>[r]_{Z}
 \ar@<0.4ex>[r]^{Z}
  &{2}
\save !<2pt,0pt>
\ar@(ur,dr)@{->}^{\ma{\bar{\lambda}^{\rr}}}
\restore}$ from $\ind'_1(\rr)$ defines the
pairs   $(\lambda ^{\rr},Z)$ and
$(\lambda ^{\rr},-Z)$ of the form \eqref{re2}; they are isomorphic via
$S=\matt{0&-1\\1&0}$ (see \eqref{kix}). Each system $M_1(\mu^{\rr})\in\ind'_1(\rr)$ defines the
pairs $(\mu^{\rr},I_2)$ and $(\mu
^{\rr},-I_2)$; they are not isomorphic since
$I_2$ and $-I_2$ are not congruent over
$\rr$.

Lemma \ref{ll} and \eqref{bxo} ensure (b$_1$) (we do not write $\nu \dd\bar{\nu}^r$ since $\nu$ is determined up to replacement by $\bar\nu$).

\item[(b$_2$):]\emph{$\varepsilon =-1$.} Let $\lambda=a+bi$ $(a,b\in\rr,\ b\ne 0)$.  If $M_{-1}(\lambda^{\rr})$ is
isomorphic to a selfdual system, then
$\lambda^r=\lambda$ or
$\lambda^r=\bar{\lambda}$. If
$\lambda^r=\lambda$, then $M_{-1}(\lambda^{\rr})$ is not isomorphic to a selfdual system; otherwise by Lemma \ref{lll} there is a nonsingular $P$ such that
 \[ \xymatrix{
 {2}  \ar[d]_{I_2}
\save !<-2pt,0pt>
 \ar@(ul,dl)@{<-}_{\matt{a&-b\\b&a}}
\restore \ar@<-0.4ex>[rr]_{-I_2}
 \ar@<0.4ex>[rr]^{I_2}
  &&{2}     \ar[d]^{P}
\save !<2pt,0pt>
\ar@(ur,dr)@{->}^{\matt{a&-b\\b&a}}
\restore
   \\
 {2}
\save !<-2pt,0pt>
 \ar@(ul,dl)@{<-}_{\matt{a&-b\\b&a}}
\restore \ar@<-0.4ex>[rr]_{-P=P^T}
 \ar@<0.4ex>[rr]^{P}
  &&{2}
\save !<2pt,0pt>
\ar@(ur,dr)@{->}^{P\matt{a&-b\\b&a}P^{-1}
=\matt{a&b\\-b&a}} \restore
 }
\] Then $P=\matt{0&-x\\x&0}$ for
some $x\ne 0$ since $P^T=-P$. The equality
$P\matt{a&-b\\b&a}=\matt{a&b\\-b&a}P$
implies that $b=0$, which contradicts our assumption that $\lambda
\notin \rr$.

If $\lambda^r=\bar{\lambda}$, then
$\matt{a&b\\-b&a}L =L\matt{a&b\\-b&a}$ with
$L:=\matt{0&-1\\1&0}$, and so
$M_{-1}(\lambda^{\rr})$ is isomorphic to
a selfdual system: \[ \xymatrix{
 {2}  \ar[d]_{I_2}
\save !<-2pt,0pt>
 \ar@(ul,dl)@{<-}_{\matt{a&-b\\b&a}}
\restore \ar@<-0.4ex>[rr]_{-I_2}
 \ar@<0.4ex>[rr]^{I_2}
  &&{2}     \ar[d]^{L}
\save !<2pt,0pt>
\ar@(ur,dr)@{->}^{\matt{a&b\\-b&a}}
\restore
   \\
 {2}
\save !<-2pt,0pt>
 \ar@(ul,dl)@{<-}_{\matt{a&-b\\b&a}}
\restore \ar@<-0.4ex>[rr]_{L^T}
 \ar@<0.4ex>[rr]^{L}
  &&{2}
\save !<2pt,0pt>
\ar@(ur,dr)@{->}^{\matt{a&b\\-b&a}}
\restore
 }
\]

Using \eqref{vfo}, we obtain \begin{align*}\label{yua}
\ind'_{-1}(\rr)=&\big\{\xymatrix{
 {2}
\save !<-2pt,0pt>
 \ar@(ul,dl)@{<-}_{\ma{\lambda^{\rr}}}
\restore \ar@<-0.4ex>[r]_{L^T}
 \ar@<0.4ex>[r]^{L}
  &{2}
\save !<2pt,0pt>
\ar@(ur,dr)@{->}^{\ma{\bar{\lambda}^{\rr}}}
\restore}
  \big|\,\lambda \in\cc^{\ud}\setminus\rr,\
\lambda ^{r}=\bar\lambda
 \big\},
                              \\
\ind''_{-1}(\rr)=& \{M_{-1}(0),\
M_{-1}(1)\}\cup
\{M_{1}(-1)\,|\,
\text{if $r$ is odd}\}
                             \\
&\cup
\{M_{-1}(\mu^{\rr})\,|\, \mu
\in\cc^{\ud}\setminus\rr,\ \mu^{r^2}=\mu
,\ \mu^{r}\ne \bar\mu\}.
\end{align*}
Each system from
$\ind'_{-1}(\rr)$ defines
the pairs $(\lambda
^{\rr},L)$ and $(\lambda ^{\rr},-L)$; they are isomorphic via
$S=\matt{0&-1\\1&0}$ (see \eqref{kix} and \eqref{re2}).

Lemma \ref{ll} and \eqref{bxo} ensure (b$_2$).
\end{itemize}

\subsection{Case {\rm(C)}: $\ff=\hh$.}

In this case,
$\varepsilon =\pm 1$, $a\mapsto \wt a$ is quaternion
conjugation
\eqref{yyb} or quaternion
semiconjugation  \eqref{ybb}, and $\C(\hh)$
is given in Lemma \ref{dcj}(c). By \eqref{tct} and
\eqref{oof},
\[\ind_{\varepsilon}(\hh)=\{M_{\varepsilon
}(\lambda)\,|\,\lambda \in\cc^{\ud},\
\lambda^{r^2}=\lambda \}.\]

If $M_{\varepsilon }(\lambda ),
M_{\varepsilon }(\mu )\in
\ind_{\varepsilon}(\hh)$, then
\begin{equation}\label{pup}
M_{\varepsilon
}(\mu)\simeq
M_{\varepsilon}(\lambda )^{\circ}\qquad
\Longleftrightarrow\qquad \mu =\lambda
^r\text{ or }\mu =\bar\lambda ^r.
\end{equation}
Indeed, if $M_{\varepsilon }(\mu)\simeq
M_{\varepsilon}(\lambda )^{\circ}$, then the $1\times 1$ matrix
$[\mu]$ is similar to
$[\lambda^r]^{\st}=[\bar\lambda ^r]$. By
Lemma \ref{dcj}(c), $\mu =\lambda^r$ or
$\mu=\bar\lambda^r$. Conversely, let
$\mu =\lambda^r$ or $\mu
=\bar\lambda^r$. We can take $\mu=\bar\lambda^r$ since $\lambda $ is determined up to
replacement by $\bar\lambda $. Then
$M_{\varepsilon }(\mu)\simeq
M_{\varepsilon}(\lambda )^{\circ}$ since
\[ \xymatrix{
 {1}  \ar[d]_{1}
\save !<-2pt,0pt>
 \ar@(ul,dl)@{<-}_{\mu}
\restore \ar@<-0.4ex>[rr]_{\varepsilon}
 \ar@<0.4ex>[rr]^{1}
  &&{1}     \ar[d]^{\varepsilon }
\save !<2pt,0pt> \ar@(ur,dr)@{->}^{\mu
^r} \restore
   \\
 {1}
\save !<-2pt,0pt>
 \ar@(ul,dl)@{<-}_{\overline{\lambda } ^r}
\restore \ar@<-0.4ex>[rr]_{1}
 \ar@<0.4ex>[rr]^{\varepsilon }
  &&{1}
\save !<2pt,0pt> \ar@(ur,dr)@{->}^{\overline{\lambda }}
\restore
 }
\]

 If $M_{\varepsilon }(\lambda )\in
\ind_{\varepsilon}(\hh)$ is
isomorphic to a selfdual system, then  by \eqref{pup} $\lambda ^r =\lambda$ or
$\lambda ^r =\bar\lambda$. Conversely,
\begin{itemize}
  \item
if $\lambda ^r =\bar\lambda$,
then $M_{\varepsilon }(\lambda )$ is
isomorphic to a selfdual system:
\begin{equation}\label{ast}
\begin{split}
\xymatrix{
 {1}  \ar[d]_{1}
\save !<-2pt,0pt>
 \ar@(ul,dl)@{<-}_{\lambda }
\restore \ar@<-0.4ex>[rr]_{\varepsilon}
 \ar@<0.4ex>[rr]^{1}
  &&{1}     \ar[d]^{\varepsilon
  \qquad\qquad\qquad
  \ma{\delta_{\varepsilon}:=}
 \begin{cases}
 1&\text{ if }\varepsilon =1,\\
 i&\text{ if }\varepsilon =-1;
 \end{cases}
  }
\save !<2pt,0pt> \ar@(ur,dr)@{->}^{\lambda
^r=\overline{\lambda}} \restore
   \\
 {1}
\save !<-2pt,0pt>
 \ar@(ul,dl)@{<-}_{\lambda }
\restore
\ar@<-0.4ex>[rr]_{\overline{\delta}_{\varepsilon}}
 \ar@<0.4ex>[rr]^{\delta_{\varepsilon} }
  &&{1}
\save !<2pt,0pt>
\ar@(ur,dr)@{->}^{\overline{\lambda}}
\restore
 }
\end{split}
\end{equation}

  \item
if $\lambda ^r =\lambda$ and $\lambda \notin\rr$ (the case $\lambda \in\rr$ is considered in the previous paragraph), then $M_{\varepsilon }(\lambda )\in
\ind_{\varepsilon}(\hh)$ is isomorphic to a selfdual system if and only if either $\varepsilon =1$ and the involution is \eqref{ybb}, or $\varepsilon =-1$ and the involution is \eqref{yyb}. Indeed, suppose that $M_{\varepsilon }(\lambda )$ is
isomorphic to a selfdual system. By Lemma \ref{lll}, there exists $h\in\hh$ such that
\begin{equation}\label{yuc}
 \begin{split}
\xymatrix{
 {1}  \ar[d]_{1}
\save !<-2pt,0pt>
 \ar@(ul,dl)@{<-}_{\lambda }
\restore \ar@<-0.4ex>[rr]_{\varepsilon}
 \ar@<0.4ex>[rr]^{1}
  &&{1}     \ar[d]^{h}
\save !<2pt,0pt> \ar@(ur,dr)@{->}^{\lambda
^r=\lambda } \restore
   \\
 {1}
\save !<-2pt,0pt>
 \ar@(ul,dl)@{<-}_{\lambda }
\restore \ar@<-0.4ex>[rr]_{ \wt{h}}
 \ar@<0.4ex>[rr]^{h}
  &&{1}
\save !<2pt,0pt> \ar@(ur,dr)@{->}^{\overline{\lambda} }
\restore
 } \end{split}
\end{equation}
is an isomorphism.
If either $\varepsilon =1$ and the involution is \eqref{ybb}, or $\varepsilon =-1$ and the involution is \eqref{yyb}, then \eqref{yuc} holds for $h=j$.
If $\varepsilon =1$ and the involution is \eqref{yyb}, then \eqref{yuc} implies $h=\bar h$, $h\in\rr$, $h\lambda =\bar \lambda h$, and so $\lambda \in\rr$, a contradiction.
If $\varepsilon =-1$ and the involution is \eqref{ybb}, then  \eqref{yuc} implies $-h=\widehat h$, $h\in\rr i$, $h\lambda =\bar \lambda h$, and so $\lambda \in\rr$, a contradiction.
\end{itemize}

\noindent
The following cases are possible:

\begin{itemize}
  \item[(c$_1$):] \emph{$\varepsilon =1$ and the involution is  quaternion
conjugation \eqref{yyb}.} Then
\begin{align*}\label{yua}
\ind'_{1}(\hh)&=
\big\{\xymatrix{
 {1}
\save !<-2pt,0pt>
 \ar@(ul,dl)@{<-}_{\lambda }
\restore
\ar@<-0.4ex>[r]_{1}
 \ar@<0.4ex>[r]^{1}
  &{1}
\save !<2pt,0pt>
\ar@(ur,dr)@{->}^{\overline{\lambda}}
\restore}
 \, \big|\,\lambda \in\cc^{\ud},\
\lambda^{r}=\bar\lambda
 \big\},
                                   \\
\ind''_{1}(\hh)&=
\{M_{1}(\mu)\,|\, \mu
\in\cc^{\ud},\ \mu ^{r^2}=\mu ,\
\mu^{r}\ne \bar\mu\}. \end{align*}

Each system from $\ind'_{1}(\hh)$ defines the pairs $(\lambda
,1)$ and $(\lambda ,-1)$; they are not isomorphic since $\bar c\,1
c\ne -1$ for all  $c\in\hh$.

Lemma \ref{ll} and \eqref{fdep} ensure (c$_1$).

 \item[(c$_2$):]  \emph{$\varepsilon =1$ and the involution is quaternion
semiconjugation \eqref{ybb}.} Then
\begin{align*}
\ind'_{1}(\hh)=&
\big\{\xymatrix{
 {1}
\save !<-2pt,0pt>
 \ar@(ul,dl)@{<-}_{\lambda }
\restore
\ar@<-0.4ex>[r]_{1}
 \ar@<0.4ex>[r]^{1}
  &{1}
\save !<2pt,0pt>
\ar@(ur,dr)@{->}^{\overline{\lambda}}
\restore}
\,
  \big|\,\lambda \in\cc^{\ud},\
\lambda^{r}=\bar\lambda
 \big\}
\\&\cup
\big\{\xymatrix{
 {1}
\save !<-2pt,0pt>
 \ar@(ul,dl)@{<-}_{\mu}
\restore
\ar@<-0.4ex>[r]_{j}
 \ar@<0.4ex>[r]^{j}
  &{1}
\save !<2pt,0pt>
\ar@(ur,dr)@{->}^{\overline{\mu}}
\restore}
\,  \big|\,\mu \in\cc^{\ud}\setminus\rr,\
{\mu}^{r}=\mu
 \big\},
                                   \\
\ind''_{1}(\hh)=&
\{M_{1}(\nu)\,|\, \nu
\in\cc^{\ud},\ \nu^{r^2}=\nu ,\ {\nu}^{r}\ne \nu,\ {\nu}^{r}\ne \bar\nu\}.
\end{align*}

Each system $\xymatrix{
 {1}
\save !<-2pt,0pt>
 \ar@(ul,dl)@{<-}_{\lambda }
\restore
\ar@<-0.4ex>[r]_{1}
 \ar@<0.4ex>[r]^{1}
  &{1}
\save !<2pt,0pt>
\ar@(ur,dr)@{->}^{\overline{\lambda}}
\restore}$ from
$\ind'_{1}(\hh)$ defines the pairs $(\lambda
,1)$ and $(\lambda ,-1)$.
\begin{itemize}
  \item
If $\lambda \notin\rr$, then  $(\lambda
,1)$ and $(\lambda ,-1)$
are not isomorphic. On the contrary, suppose that there is a nonzero $c\in\hh$ such that
\begin{equation}\label{wwmy}
c^{-1}\lambda c=\lambda,\qquad \widehat c\,1 c=-1
\end{equation}
(see \eqref{kix}). By $c^{-1}\lambda c=\lambda$, we have $c\in\cc$, which contradicts $\widehat c\,1c=-1$.

  \item If $\lambda \in\rr$, then
$(\lambda ,1)$ and $(\lambda ,-1)$ are isomorphic since \eqref{wwmy} holds for $c=j$.
\end{itemize}

The pairs $(\mu
,j)$ and $(\mu ,-j)$ constructed by
$\xymatrix{
 {1}
\save !<-2pt,0pt>
 \ar@(ul,dl)@{<-}_{\mu}
\restore
\ar@<-0.4ex>[r]_{j}
 \ar@<0.4ex>[r]^{j}
  &{1}
\save !<2pt,0pt>
\ar@(ur,dr)@{->}^{\overline{\mu}}
\restore}$  from
$\ind'_{1}(\hh)$
are  isomorphic via $i$ since $i\mu =\mu i$ and $\widehat \imath\,ji=-iji=-ki=-j$.

Lemma \ref{ll} and
\eqref{fdep} ensure (c$_2$).

 \item[(c$_3$):]  \emph{$\varepsilon =-1$ and the involution is quaternion
conjugation \eqref{yyb}.} Then
\begin{align*}
\ind'_{-1}(\hh)=&
\big\{\xymatrix{
 {1}
\save !<-2pt,0pt>
 \ar@(ul,dl)@{<-}_{\lambda }
\restore
\ar@<-0.4ex>[r]_{-i}
 \ar@<0.4ex>[r]^{i}
  &{1}
\save !<2pt,0pt>
\ar@(ur,dr)@{->}^{\overline{\lambda}}
\restore}
\,  \big|\,\lambda \in\cc^{\ud},\
\lambda^{r}=\bar\lambda
 \big\}
\\&\cup
\big\{\xymatrix{
 {1}
\save !<-2pt,0pt>
 \ar@(ul,dl)@{<-}_{\mu}
\restore
\ar@<-0.4ex>[r]_{-j}
 \ar@<0.4ex>[r]^{j}
  &{1}
\save !<2pt,0pt>
\ar@(ur,dr)@{->}^{\overline{\mu}}
\restore}\,
  \big|\,\mu \in\cc^{\ud}\setminus\rr,\
{\mu}^{r}=\mu
 \big\},
                                   \\
\ind''_{-1}(\hh)=&
\{M_{-1}(\nu)\,|\, \nu
\in\cc^{\ud},\ \nu ^{r^2}=\nu ,\
{\nu}^{r}\ne \nu,\
{\nu}^{r}\ne \bar\nu\}.
\end{align*}

If $\lambda \notin\rr$,  then the pairs
$(\lambda ,i)$ and $(\lambda ,-i)$ constructed by a system from
$\ind'_{-1}(\hh)$ are not isomorphic. On the contrary, suppose there exists  $c\in\hh$ such that
\begin{equation}\label{wwm}
c^{-1}\lambda c=\lambda,\qquad \bar c\,i c=-i
\end{equation}
(see \eqref{kix}). Since $\lambda =\lambda _1+\lambda _2i$ with $\lambda _1,\lambda _2\in\rr$ and $\lambda _2\ne 0$, the equality $c^{-1}\lambda c=\lambda $ implies that $ic=ci$, and so $c\in\cc$, which contradicts $\bar c\,ic=-i$.
  If $\lambda \in\rr$, then
$(\lambda ,i)$ and $(\lambda ,-i)$ are isomorphic since \eqref{wwm} holds for $c=j$.

The pairs $(\mu
,j)$ and $(\mu ,-j)$ are  isomorphic via $i$ since $i\mu =\mu i$ and $\bar \imath\,ji=-iji=-ki=-j$.

Lemma \ref{ll} and
\eqref{fdep} ensure (c$_3$).

 \item[(c$_4$):]  \emph{$\varepsilon =-1$ and the involution is quaternion
semiconjugation \eqref{ybb}.} Then
\begin{align*}
\ind'_{-1}(\hh)&=
\big\{\xymatrix{
 {1}
\save !<-2pt,0pt>
 \ar@(ul,dl)@{<-}_{\lambda }
\restore
\ar@<-0.4ex>[r]_{-i}
 \ar@<0.4ex>[r]^{i}
  &{1}
\save !<2pt,0pt>
\ar@(ur,dr)@{->}^{\overline{\lambda}}
\restore}
\,  \big|\,\lambda \in\cc^{\ud},\
\lambda^{r}=\bar\lambda\big\},
                                   \\
\ind''_{-1}(\hh)&=
\{M_{-1}(\mu)\,|\, \mu
\in\cc^{\ud},\ \mu ^{r^2}=\mu ,\
\mu^{r}\ne \bar\mu\}. \end{align*}

The pairs $(\lambda ,i)$ and
$(\lambda ,-i)$ constructed by a system from $\ind'_{-1}(\hh)$
are not isomorphic since $\widehat c\,ic=i\bar c c\ne -i$ for all
$c\in\hh$.

Lemma \ref{ll} and
\eqref{fdep} ensure (c$_4$).
\end{itemize}

\section*{Acknowledgements}
This paper is a result of a student seminar held at the University of S\~ao Paulo during the visit of V.V.~Sergeichuk in 2019 and 2020;
he is grateful to the university for hospitality and to the FAPESP
for financial support (2018/24089-4).

\end{document}